\documentstyle[11pt]{article}
\def\a{\alpha}            \def\b{\beta}
\def\g{\gamma}            
\def\d{\delta}            \def\D{\Delta}
\def\th{\theta}           
       \def\ep{\epsilon}
\def\vep{\varepsilon}

\def\sg{\sigma}           \def\Sg{\Sigma}
              
          \def\Lm{\Lambda}

\def\h{\hat}              \def\wh{\widehat}
\def\tl{\tilde}

                  \def\cF{{\cal F}}

\def\lb{\label}                 \def\ot{\otimes}
\def\cd{\cdots}                 \def\ld{\ldots}

\def\fs{\footnotesize}          
            \def\ls{\large}

\def\r#1{(\ref{#1})}

\def\ZZ{Z\!\!\!\!\!\!Z}

\def\Z2{Z\!\!\!\!\!\!Z_{\,2}}
\def\CC{{\bf C}\!\!\!{\rm l}\,}

\newcommand{\bn}{\begin{equation}}
\newcommand{\ed}{\end{equation}}
\newcommand{\bneqn}{\begin{eqnarray}}
\newcommand{\edeqn}{\end{eqnarray}}
\newcommand{\bnth}{\begin{theorem}}
\newcommand{\edth}{\end{theorem}}
\newcommand{\bnpr}{\begin{proposition}}
\newcommand{\edpr}{\end{proposition}}
\newcommand{\bnlm}{\begin{lemma}}
\newcommand{\edlm}{\end{lemma}}
\newcommand{\bndef}{\begin{definition}}
\newcommand{\eddef}{\end{definition}}
\newcommand{\bntb}{\begin{tabular}}
\newcommand{\edtb}{\end{tabular}}
\def\nin{\noindent}
\newtheorem{definition}{Definition}[section]
\newtheorem{lemma}{Lemma}[section]
\newtheorem{proposition}{Proposition}[section]
\newtheorem{theorem}{Theorem}[section]

\begin{document}

\begin{center}
\vspace*{0.8cm}

{\LARGE{\bf Modified Affine Hecke Algebras and\\
Drinfeldians of Type A}
\footnote{The talk given by V.N. Tolstoy}
}

\vskip 1.3cm
{\large{\bf V.N. Tolstoy$^{1}$, O.V. Ogievetsky$^{2}$,\\[3pt]
P.N. Pyatov$^{3}$ {\rm and} A.P. Isaev$^{3}$}}

\vskip 0.5cm
{\it $^{1}$Institute of Nuclear Physics,
Moscow State University\\
119899 Moscow \& Russia\\
(e-mail: tolstoy@anna19.npi.msu.su)
}

\vskip 0.15cm
{\it $^{2}$Marseille University and Center of Theoretical Physics, CNRS\\
Luminy - Case 907-13288 Marseille Cedex 9 \& France\\
(e-mail: oleg@cptsu5.univ-mrs.fr)
}

\vskip 0.15cm
{\it $^{3}$Bogoliubov Laboratory of Theoretical Physics, JINR\\
Joint Institute of Nuclear Reserch\\
141980 Dubna, Moscow region \& Russia\\
(e-mails: pyatov@thsun1.jinr.ru, isaevap@thsun1.jinr.ru)
}

\end{center}

\vspace{0,7 cm}
\begin{abstract}
We introduce a modified affine Hecke algebra by a singular
transformation of the usual affine Hecke algebra $\h{H}_{q}(l)$ of
type $A_{l-1}$. The modified affine Hecke algebra $\h{H}_{q\eta}({l})$
($\h{H}^{+}_{q\eta}({l})$) depends on two deformation parameters $q$
and $\eta$. When the parameter $\eta$ is equal to zero the algebra
$\h{H}_{\!q\eta=0}(l)$ coincides with $\h{H}_{q}(l)$, if the parameter
q goes to 1 the algebra $\h{H}^{+}_{\!q=\!1\eta}(l)$ is isomorphic to
the degenerate affine Hecke algebra $\Lm_{\eta}(l)$ introduced by
Drinfeld. We construct a functor ${\cal F}_{\!q\eta}$ from a category of
representations of $H_{q\eta}^{+}(l)$ into a category of representations
of Drinfeldian $D_{q\eta}(sl(n\!+\!1))$ which has been introduced
by the first author. This functor depends on two continuous
deformation parameters $q$ and $\eta$. If the parameter $\eta$ is
equal to zero then the functor ${\cal F}_{\!q\eta=0}$ coincides with
the duality functor constructed by Chari and Pressley for the affine
Hecke algebra $\h{H}_{q}^{+}(l)$ and the quantum affine algebra
$U_{q}(sl(n\!+\!1)[u])$. When the parameter $q$ goes to 1 the
functor ${\cal F}_{\!q=\!1\eta}$ coincides with Drinfeld's
functor for the degenerate affine Hecke algebra $\Lm_{\eta}(l)$
and the Yangian $Y_{\eta}(sl(n\!+\!1))$.
\end{abstract}

\setcounter{equation}{0}
\section{Introduction}

One of the most remarkable results of the classical representation
theory is the Frobenius-Schur duality between the finite-dimensional
irreducible representations of the general or special linear groups
and symmetric groups. The duality means that
any finite-dimensional irreducible representation of the Lie algebra
$g$ (or its universal enveloping algebra $U(g)$),
where $g=gl(n\!+\!1)$ or $sl(n\!+\!1)\simeq A_n$,
can be obtained by decomposing of the $l$-fold tensor product of the
fundamental (natural) representation $V=\CC^{n+1}$ with respect to
the action of the symmetric group $S(l)$ (or its group algebra
$\CC[S(l)]$).

After discovery of the quantum groups \cite{D1,J1}, Jimbo
\cite{J2} proved the q-analog of the Frobenius- Schur duality
replacing $U(g)$ by $U_q(g)$ and $\CC[S(l)]$ by its q-analogue
$H_q(l)$, the Hecke algebra of type $A_{l-1}$. Slightly earlier
in 1985, Drinfeld \cite{D2} discovered an analogue of the
Frobenius-Schur theory for the Yangian $Y_{\eta}(sl(n\!+\!1))$
and the degenerate affine Hecke algebra $\Lm_{\eta}(l)$.
Later, Chari and Pressley \cite{ChP1}
proved the q-analogue of the duality for the quantum affine
algebra $U_q(\wh{sl}(n\!+\!1))$ and the affine Hecke algebra
$\h{H}_q(l)$.

In this paper, we extend the results of Drinfeld and Chari-Pressley
to the case of the Drinfeldian $D_{\!q\eta}(sl(n\!+\!1))$
\cite{T2} - \cite{T4} which is the rational-trigonometric deformation
of the universal enveloping algebra of the loop algebra
$sl(n\!+\!1)[u]$. In this case, the role of $\h{H}_q(l)$ is played
the modified affine Hecke algebra $\h{H}_{\!q\eta}^{+}(l)$ which we
obtain by a singular transformation of the affine Hecke $\h{H}_q(l)$.
Our functor ${\cal F}_{\!q\eta}$ from a category of representations of
$H_{\!q\eta}^{+}(l)$ in a category of those of the Drinfeldian
$D_{\!q\eta}(sl(n\!+\!1))$ depends on two continuous
deformation parameters $q$ and $\eta$. If the parameter $\eta$ is
equal to zero then the functor ${\cal F}_{\!q\eta=0}$ coincides with
the duality functor constructed by Chari and Pressley \cite{ChP1}
for the affine Hecke algebra $\h{H}_{q}^{+}(l)$ and the quantum
affine algebra $U_{q}(sl(n\!+\!1)[u])$. When the parameter $q$
goes to 1 the functor ${\cal F}_{\!q=\!1\eta}$ coincides with
Drinfeld's functor for the degenerate affine Hecke algebra
$\Lm_{\eta}(l)$ and the Yangian $Y_{\eta}(sl(n\!+\!1))$ \cite{D2}.

\setcounter{equation}{0}
\section{Affine Hecke and modified affine Hecke algebras}

We start from the definition of the affine Hecke algebra
\cite{BZ,ChP2,R}.
\begin{definition}
The affine Hecke algebra $\hat{H}_{q}(l)\!:=\!\hat{H}_{q}(A_{l-l})$
of type $A_{l-1}$ is an associative algebra over $\CC[q,q^{-1}]$, generated
by the elements $\sg_{1}^{\pm 1},\sg_{2}^{\pm 1},\ld,\sg_{l-1}^{\pm 1}$,
and $z_{1}^{\pm 1},z_{2}^{\pm 1},\ld,z_{l}^{\pm 1}$ with the following
defining relations:
\bneqn
\sg_{i}^{}\sg_{i}^{-1}&=&\sg_{i}^{-1}\sg_{i}^{}\;=\;1~,
\lb{AH1}
\\[7pt]
\sg_{i}^{}-\sg_i^{-1}&=&(q-q^{-1})~,
\lb{AH2}
\\[7pt]
\sg_{i}\sg_{i+1}\sg_{i}&=&\sg_{i+1}\sg_{i}\sg_{i+1}~,
\lb{AH3}
\\[7pt]
\sg_{i}\sg_{j}&=&\sg_{j}\sg_{i} \qquad\quad {\rm if}\;\;|i-j|>1~,
\lb{AH4}
\\[7pt]
z_{j}^{}z_{j}^{-1}&=&z_{j}^{-1}z_{j}^{}~=~1~,
\lb{AH5}
\\[7pt]
z_{j}z_{k}&=&z_{k}z_{j}~,
\lb{AH6}
\\[7pt]
\sg_{i}z_{j}&=&z_{j}\sg_{i}\qquad\quad{\rm if}\;\;j\neq i\;{\rm or}\;i+1~,
\lb{AH7}
\\[7pt]
\sg_{i}^{}z_{i}^{}&=&z_{i+1}^{}\sg_{i}^{-1}~.
\lb{AH8}
\edeqn
An associative algebra generated by the elements
$\sg_{i}^{\pm 1}$, $i\in \{1,2,\ld,l-1\}$,
with
the defining relations
\r{AH1}--\r{AH4} is called the Hecke algebra
$H_{q}(l)\!:=\!H_{q}(A_{l-1})$.
\lb{AHD1}
\end{definition}
Sometimes it is useful to use the last relation (\ref{AH8}) in another
forms. Namely applying the relation (\ref{AH2}) one obtains
\bn
\sg_{i}z_{i}-z_{i+1}\sg_{i}=(q^{-1}-q)z_{i+1}
\lb{AH9}
\ed
or
\bn
z_{i}\sg_{i}-\sg_{i}z_{i+1}=(q^{-1}-q)z_{i+1}~.
\lb{AH10}
\ed
The permutation relations for the inverse powers of the generators
$z_i$ looks like
\bn
\begin{array}{rcl}
z_{i}^{-1}\sg_{i}-\sg_{i}z_{i+1}^{-1}&=&(q^{-1}-q)z_{i}^{-1}~,
\\[7pt]
\sg_{i}z_{i}^{-1}-z_{i+1}^{-1}\sg_{i}&=&(q^{-1}-q)z_{i}^{-1}~.
\end{array}
\lb{AH11}
\ed

Using the relations (\ref{AH9})-(\ref{AH11}) and (\ref{AH7}) it is
easy to see that any polynomial of the elements $\sg_{i}^{\pm 1}$
$(i=1,\ld,l-1)$, and $z_{j}^{\pm 1}$ $(j=1,\ld,l)$
may be put in order such that all elements $\sg_{i}^{\pm 1}$ are located
from the left-hand side (or from the right-hand side) of the elements
$z_{j}^{\pm 1}$, i.e. any polynomials of $\sg_{i}^{\pm 1}$ and
$z_{j}^{\pm 1}$ is represented as a sum of the monomials of the type
\bn
z_{1}^{n_{1}}z_{2}^{n_{2}}\cd z_{l}^{n_{l}}
\sg_{i_{1}}\sg_{i_{2}}\cd\sg_{i_{k}}\qquad
({\rm or}\quad\sg_{i_{1}}\sg_{i_{2}}\cd\sg_{i_{k}}
z_{1}^{n_{1}}z_{2}^{n_{2}}\cd z_{l}^{n_{l}})~,\qquad n_i\in \ZZ\,,
\lb{AH12}
\ed
where among the elements $\sg_{i_{j}}$ can be equal.
This result is reformulated as the following proposition.
\bnpr
There is an isomophism of the vector spaces $\h{H}_{q}(l)$ and
$\CC[z_{1}^{\pm 1},\ld,z_{l}^{\pm 1}]\ot H_{q}(l)$
(or $H_{q}(l)\ot \CC[z_{1}^{\pm 1},\ld,z_{l}^{\pm 1}]$):
\bn
\h{H}_{q}(l)\simeq\CC[z_{1}^{\pm 1},\ld,z_{l}^{\pm 1}]\ot H_{q}(l)
\qquad\quad ({\rm or}\;\; \h{H}_{q}(l)\simeq H_{q}(l)\ot
\CC[z_{1}^{\pm 1},\ld,z_{l}^{\pm 1}])~.
\lb{AH13}
\ed
\lb{PAH1}
\edpr
\vskip -10pt

The subalgebra  $\h{H}_{q}^{+}(l)\subset \h{H}_{q}(l)$, which is
generated by $H_{q}(l)$ and the elements $z_{1},z_{2},\ld,z_{l}$
will be also called the affine Hecke algebra.

The affine Hecke $\h{H}_{q}^{+}(l)$ (and also $\h{H}_{q}(l)$)
does not contain any singular elements at $q\to 1$ and
\bn
\lim_{q\to 1}\h{H}_{q}(l)\simeq\h{\Sg}(l)\ , \qquad \mbox{and}\qquad
\lim_{q\to 1}\h{H}_{q}^+(l)\simeq\h{\Sg}^+(l)\ ,
\lb{AH14}
\ed
where by $\h{\Sg}(l)$ ($\h{\Sg}^+(l)$) we denote the affine symmetric
group algebra generated by the group algebra of the symmetric group
$\CC[S(l)]$ and the affine elements
$z_{1}^{\pm},z_{2}^{\pm},\ld,z_{l}^{\pm}$ ($z_{1},z_{2},\ld,z_{l}$)
with the defining relation (\ref{AH1})-(\ref{AH8}) for $q=1$.

Now we introduce a modified the affine Hecke algebra by
the singular translation of the affine elements $z_{j}$:
\bn
u_{j}=z_{j}+\frac{\eta}{q-q^{-1}}\qquad\quad {\rm for}\;j=1,2,\ld,l~.
\lb{AH15}
\ed
This transformation  changes only the last relation (\ref{AH8}))
from the set \r{AH1}--\r{AH8}, which takes now the form
\bn
\sg_{i}^{}u_{i}^{}=u_{i+1}^{}\sg_{i}^{-1}+\eta~.
\lb{AH16}
\ed
A remarkable fact is that while the transformation \r{AH15}
contains terms which are singular, in the classical limit
$q\rightarrow 1$, the permutation relations \r{AH16} for the newly
defined generators $u_i$ do not. So we have:
\begin{definition}
The modified affine Hecke algebra
$\hat{H}_{q\eta}^{+}(l)\!=\!\hat{H}_{q\eta}^{+}(A_{l-1}) $ of
type $A_{l-1}$ is an associative algebra over $\CC[q,q^{-1},\eta]$
generated by the elements
$\sg_{1}^{\pm 1},\sg_{2}^{\pm 1},\ld,\sg_{l-1}^{\pm 1}$,
and $u_{1}^{},u_{2}^{},\ld,u_{l}^{}$ with the following
defining relations:
\bneqn
\sg_{i}^{}\sg_{i}^{-1}&=&\sg_{i}^{-1}\sg_{i}^{}\;=\;1~,
\lb{AH17}
\\[7pt]
\sg_{i}^{}-\sg^{-1}&=&(q-q^{-1})~,
\lb{AH18}
\\[7pt]
\sg_{i}\sg_{i+1}\sg_{i}&=&\sg_{i+1}\sg_{i}\sg_{i+1}~,
\lb{AH19}
\\[7pt]
\sg_{i}\sg_{j}&=&\sg_{j}\sg_{i}\qquad\qquad{\rm if}\;\;|i-j|>1~,
\lb{AH20}
\\[7pt]
u_{j}u_{k}&=&u_{k}u_{j}~,
\lb{AH21}
\\[7pt]
\sg_{i}u_{j}&=&u_{j}\sg_{i}\qquad\qquad{\rm if}\;\;j\neq i\;{\rm or}\;i+1~,
\lb{AH22}
\\[7pt]
\sg_{i}^{}u_{i}^{}&=&u_{i+1}^{}\sg_{i}^{-1}+\eta~.
\lb{AH23}
\edeqn
\lb{AHD2}
\end{definition}

\vskip-15pt
The "$\eta-analog$" of the relations (\ref{AH9}), (\ref{AH10})
now looks like
\bn
\begin{array}{lcr}
\sg_{i}u_{i}-u_{i+1}\sg_{i}&=&(q^{-1}-q)u_{i+1}+\eta~,
\\[7pt]
u_{i}\sg_{i}-\sg_{i}u_{i+1}&=&(q^{-1}-q)u_{i+1}+\eta~.
\lb{AH24}
\end{array}
\ed
It is obvious  that the statement of the Proposition \ref{PAH1}
remains valid for the modified affine Hecke algebra.

One can extend the algebra $\hat{H}_{q\eta}^{+}(l)$
adding generators $u_{j}^{-1}$ inverse to the elements $u_{j}^{}$:
$u_{j}^{}u_{j}^{-1}=u_{j}^{}u_{j}^{-1}=1$. In this way one obtains
the total modified affine Hecke algebra $\hat{H}_{q\eta}(l)$.
However in the present paper we need only the subalgebra
$\hat{H}_{q\eta}^{+}(l)\subset\hat{H}_{q\eta}(l)$.

The algebra $\hat{H}_{q\eta}^{+}(l)$ is a two-parameter
$(q,\eta)$-deformation of $\h{\Sg}^+(l)$. However it is easy to see
that the modified affine Hecke algebra $\hat{H}_{q\eta}^{+}(l)$
is essentially independent of the parameter $\eta$, provided that
$\eta\neq 0$. In fact, if $\eta\neq 0$ and $\eta'\neq 0$ the map
$\h{H}_{q\eta}^{+}(l)\to \h{H}_{q\eta'}^{+}(l)$ given by
$\sg_{i}\mapsto\sg_{i}$, $\eta^{-1} u_{j}\mapsto{\eta'}^{-1}u_{j}$
is clearly an isomorphism of these algebras. Thus one might as well
take $\eta=1$, however we keep the parameter $\eta$ for visualization.

It is obvious that
$\h{H}_{\!q\eta=0}^{+}(l)=\h{H}_{q}^{+}(l)$.
On the other hand, in the limit
$q\rightarrow 1$ the modified affine Hecke algebra
goes into the degenerate affine Hecke
algebra $\Lm_{\eta}(l)$ constructed by Drinfeld in 1985 \cite{D2}
\footnote{\nin This algebra was also obtained by Drinfeld from the
affine Hecke algebra $\h{H}^{+}_{q}(l)$ by letting $q\to 1$ in
a certain non-trivial fashion.}.
The relations between the modified affine Hecke algebra
$\h{H}_{q\eta}^{+}(l)$ and the algebras $\h{H}_{q}^{+}(l)$,
$\Lm_{\eta}(l)$, $\h{\Sg}^{+}(l)$ (and also their subalgebras)
are shown in the picture:
\bn
\mbox{\begin{picture}(100,60)
\put(00,25){$\h{H}_{q\eta}^{+}(l)$}
\put(00,-25){\mbox{$\Lm_{\eta}(l)$}}
\put(80,25){\mbox{$\h{H}_{q}^{+}(l)$}}
\put(80,-25){\mbox{$\h{\Sg}^{+}(l)$}}
\put(35,29){\vector(1,0){40}}
\put(35,-21){\vector(1,0){40}}
\put(4,16){\vector(0,-1){25}}
\put(85,16){\vector(0,-1){25}}
\put(-40,25){\mbox{$H_{q}(l)\subset$}}
\put(-40,-25){\mbox{$\Sg\,(l)\,\subset$}}
\put(112,25){\mbox{$\supset H_{q}(l)$}}
\put(110,-25){\mbox{$\supset\Sg\,(l)$}\,\,.}
\put(-20,0){\mbox{\fs$q\!\to\!1$}}
\put(88,0){\mbox{\fs$q\!\to\!1$}}
\put(40,33){\mbox{\fs$\eta\!\to\!0$}}
\put(40,-17){\mbox{\fs$\eta\!\to\!0$}}
\end{picture}}
\lb{AH25}
\ed
\\[2pt]
\begin{quote}
{\fs Fig.1. A diagram of the limit algebras
of the modified affine Hecke algebra $\h{H}_{q\eta}^{+}(l)$\\
\phantom{Fig.1.} and their subalgebras. The arrows show passages
to the limits.}
\end{quote}
\vskip 10pt

\setcounter{equation}{0}
\section{Drinfeldian and Yangian of type $A_{n}$}

First we recall the defining relations of the q-quantized universal
enveloping algebra $U_q(sl(n\!+\!1))$
($sl(n\!+\!1):=sl(n\!+\!1,\CC)\simeq A_{n}$)
and construction of its Cartan-Weyl basis.

Let $\Pi:=\{\a_{1},\ld,\a_{n}\}$ be a system of simple roots of
$sl(n+1)$ endowed with the following scalar product:
$(\a_{i},\a_{j})=(\a_{j},\a_{i})$, $(\a_{i},\a_{i})=2$,
$(\a_{i},\a_{i+1})=-1$, $(\a_i,\a_j)=0$ $((|i-j|>1)$.
The corresponding Dynkin diagram is presented on the picture:
\\[0pt]
\bn
\mbox{\begin{picture}(120,10)
\put(-40,0){\circle{5}}
\put(-37,0){\line(1,0){44}}
\put(-43,-15){\fs$\a_{1}$}
\put(10,0){\circle{5}}
\put(13,0){\line(1,0){35}}
\put(7,-15){\fs$\a_{2}$}
\put(55,0){$\ld$}
\put(107,0){\line(-1,0){35}}
\put(110,0){\circle{5}}
\put(113,0){\line(1,0){44}}
\put(100,-15){\fs$\a_{n-1}$}
\put(160,0){\circle{5}}
\put(157,-15){\fs$\a_{n}$}
\end{picture}}
\lb{DY1}
\ed\\[2pt]
\centerline{\fs Fig.3. Dynkin diagram of the Lie algebra $sl(n+1)$.}
\\[7pt]

The quantum algebra
$U_{q}(sl(n\!+\!1))$ is generated by the Chevalley
elements $q^{\pm h_{\a_i}}$, $e_{\pm\a_i}$ $(i=1,2,\ld,n)$ with
the defining relations:
\bn
\begin{array}{rcl}
q^{h_{\a_i}}q^{-h_{\a_i}}&=&q^{-h_{\a_i}}q^{h_{\a_i}}=1~,
\\[7pt]
q^{h_{\a_i}}q^{h_{\a_j}}&=&q^{h_{\a_j}}q^{h_{\a_i}}~,
\\[7pt]
q^{h_{\a_i}}e_{\pm\a_j}q^{-h_{\a_i}}&=&q^{\pm(\a_i,\a_j)}e_{\pm\a_j}~,
\\[7pt]
[e_{\a_i},e_{-\a_j}]&=&\d_{ij}\,
[h_{\alpha_i}]_q
\\[12pt]
[e_{\pm\a_i},e_{\pm\a_j}]&=&0\qquad\quad(|i-j|\geq 2)~,
\\[7pt]
[[e_{\pm\a_i}e_{\pm\a_j}]_{q}^{}e_{\pm\a_j}]_{q}^{}&=&0
\qquad\quad(|i-j|=1)~,
\lb{DY2}
\end{array}
\ed
where $[h]_q:= (q^{h}-q^{-h})/(q-q^{-1})$ is standard notation for
the "q-number" and $[\,\cdot\,,\,\cdot\,]_{q}$ is the q-commutator:
\bn
[e_{\b},e_{\g}]_{q}:=e_{\b}e_{\g}-q^{(\b,\g)}e_{\g}e_{\b}~.
\lb{DY3}
\ed
The Hopf structure on $U_{q}(sl(n\!+\!1))$ is given by
the following formulas for a comultiplication $\D_{q}$,
an antipode $S_{q}$, and a co-unit $\vep_{q}$:
\bn
\begin{array}{rcl}
\D_{q}(q^{\pm h_{\a_i}})&=&q^{\pm h_{\a_i}}\ot q^{\pm h_{\a_i}} ~,
\\[5pt]
\D_{q}(e_{\a_i})&=&e_{\a_i}\ot 1+q^{-h_{\a_i}}\ot e_{\a_i}~,
\\[5pt]
\D_{q}(e_{-\a_i})&=&e_{-\a_i}\ot q^{h_{\a_i}}+1\ot e_{-\a_i}~;
\end{array}
\lb{DY4}
\ed
\\[-10pt]
\bn
\begin{array}{rcl}
S_{q}(q^{\pm h_{\a_i}})&=&q^{\mp h_{\a_i}}~,\qquad\qquad\qquad\quad
\\[5pt]
S_{q}(e_{\a_i})&=&-q^{h_{\a_i}}e_{\a_i}~,
\\[5pt]
S_{q}(e_{-\a_i})&=&-e_{-\a_i}q^{-h_{\a_i}}~;
\end{array}
\lb{DY5}
\ed
\\[-10pt]
\bn
\begin{array}{rcl}
\vep_{q}(q^{\pm h_{\a_i}})&=&1~,\qquad\qquad\qquad\qquad\;
\\[5pt]
\vep_{q}(e_{\pm\a_i})&=&0~.
\end{array}
\lb{DY6}
\ed

Below we shall also use another basis in the Cartan subalgebra
of the Lie algebra $sl(n\!+\!1)$. Namely we set
\bn
\begin{array}{rcl}
e_{11}^{}&=&\mbox{\ls$\frac{1}{n+1}$}\left(nh_{\a_1}+(n\!-\!1)h_{\a_2}+\cd+
2h_{\a_{n-1}}+h_{\a_{n}}+N\right)~,
\\[7pt]
e_{22}^{}&=&\mbox{\ls$\frac{1}{n+1}$}\left(nh_{\a_1}+(n\!-\!1)h_{\a_2}+
\cd+2h_{\a_{n-1}}+h_{\a_{n}}+N\right)-h_{\a_1}~,
\\[3pt]
\ld&\ld&\ld\ld\ld\ld\ld\ld\ld\ld\ld\ld\ld\ld\ld\ld\ld\ld\ld\ld
\\[3pt]
e_{ii}^{}&=&\mbox{\ls$\frac{1}{n+1}$}\left(nh_{\a_1}+(n\!-\!1)h_{\a_2}+
\cd+2h_{\a_{n-1}}+h_{\a_{n}}+N\right)-\sum\limits_{k=1}^{i-1}h_{\a_k}~,
\\[3pt]
\ld&\ld&\ld\ld\ld\ld\ld\ld\ld\ld\ld\ld\ld\ld\ld\ld\ld\ld\ld\ld
\\[3pt]
e_{n\!+\!1n\!+\!1}^{}&=&\mbox{\ls$\frac{1}{n+1}$}\left(-h_{\a_1}-
2h_{\a_2}-\cd-
(n\!-\!1)h_{\a_{n-1}}-nh_{\a_{n}}+N\right)~.
\end{array}
\lb{DY7}
\ed
Here $N$ is a central element of $g$ (and also of $U_{q}(g)$),
which is equal to 0 for the case $g=sl(n\!+\!1)$ and $N\neq 0$ for
$g=gl(n\!+\!1)$. It is easy to see that
\bn
\begin{array}{rcl}
h_{\a_i}&=&e_{ii}-e_{i\!+\!1i\!+\!1}\qquad(i=1,\ld,n)~,
\\[7pt]
N&=&e_{11}+e_{22}+\ld+e_{n\!+\!1n\!+\!1}~.
\end{array}
\lb{DY8}
\ed
A dual basis to the elements $e_{ii}$ ($i=1,2,\ld, n\!+\!1$)
will be denoted by $\ep_i$ ($i=1,2,\ld, n\!+\!1$):
$\ep_i(e_{jj})=(\ep_i,\ep_j)=\d_{ij}$. In the terms of $\ep_i$
the positive root system $\D_{+}$ of $sl(n\!+\!1)$
is presented as follows
\bn
\D_{+}=\{\ep_i-\ep_j\,|\,1\le i<j\le n+1\}~,
\lb{DY9}
\ed
where $\ep_{i}-\ep_{i+1}$ are the simple roots:
\bn
\a_{i}=\ep_{i}-\ep_{i+1}\quad (i=1,2,\ld,n)~.
\lb{DY10}
\ed
The root $\th:=\ep_{1}-\ep_{n+1}$ is maximal one:
\bn
\th=\a_1+\a_2+\ld+\a_{n}~.
\lb{DY11}
\ed
For the root vectors $e_{\ep_i-\ep_j}$ $(i\neq j)$ the standard
notations are also used
\bn
e_{ij}:=e_{\ep_i-\ep_j}~,\qquad e_{ji}:=e_{\ep_j-\ep_i}
\qquad (1\le i<j\le n+1)~.
\lb{DY12}
\ed
In particular, $e_{ii\!+\!1}$, $e_{i\!+\!1i}$ are the Chevalley elements:
$e_{ii\!+\!1}=e_{\a_i}$, $e_{i\!+\!1i}=e_{-\a_i}$ ($i=1,\ld, n$).

For construction of the composite root vectors $e_{ij}$ $(j\neq i\pm 1)$
we fix the following normal ordering of the positive root system
$\D_{+}$ (see \cite{T1,KT1})
\bn
(\ep_1-\ep_2),(\ep_1-\ep_3,\ep_2-\ep_3),\ld,(\ep_1-\ep_i,\ld,
\ep_{i-1}-\ep_{i}),\ld,(\ep_1-\ep_{n+1},\ld,\ep_n-\ep_{n+1})\ .
\lb{DY13}
\ed
According to with this ordering we set
\bn
e_{ij}:=[e_{ik},e_{kj}]_{q^{-1}}~,\qquad
e_{ji}:=[e_{jk},e_{ki}]_{q}\qquad (1\le i<k<j\le n+1)~.
\lb{DY14}
\ed
It should be stressed that the structure of the composite root
vectors (\ref{DY14}) is independent of choice of the index $k$
in the r.h.s. of the definition \r{DY14}.
In particular one has
\bn
\begin{array}{rcccll}
e_{ij}&\!\!:=\!\!&[e_{ii\!+\!1},e_{i\!+\!1j}]_{q^{-1}}&\!\!=\!\!&
[e_{ij\!-\!1},e_{j\!-\!1j}]_{q^{-1}}\qquad &(1\le i<j\le n+1)~,
\\[7pt]
e_{ji}&\!\!:=\!\!&[e_{ji\!+\!1},e_{i\!+\!1i}]_{q}&\!\!=\!\!&
[e_{jj\!-\!1},e_{j\!-\!1i}]_{q}\qquad &(1\le i<j\le n+1)~.
\end{array}
\lb{D16}
\ed
General properties of the Cartan-Weyl basis $\{e_{ij}\}$ can be
found in \cite{T1,KT1,KT2}.

As it was noted in \cite{T2} the Dynkin diagrams of the non-twisted
affine algebras can be also used for classification of the Drinfeldians
and the Yangians. In the case of $sl(n\!+\!1)$, the Dynkin diagram of
the corresponding affine Lie algebra $\wh{sl}(n\!+\!1)$ is presented by
the picture:
\\[-7pt]

\bn
\mbox{\begin{picture}(120,60)
\put(-40,0){\circle{5}}
\put(-37,0){\line(1,0){44}}
\put(-43,-15){\fs$\a_{1}$}
\put(10,0){\circle{5}}
\put(13,0){\line(1,0){35}}
\put(7,-15){\fs$\a_{2}$}
\put(53,0){$\ld$}
\put(107,0){\line(-1,0){35}}
\put(110,0){\circle{5}}
\put(113,0){\line(1,0){44}}
\put(100,-15){\fs$\a_{n-1}$}
\put(160,0){\circle{5}}
\put(157,-15){\fs$\a_{n}$}
\put(60,49){\circle{5}}
\put(-38,2){\line(2,1){95}}
\put(158,2){\line(-2,1){95}}
\put(50,60){\fs$\d\!-\!\th$}
\end{picture}}
\lb{DY17}
\ed\\[2pt]
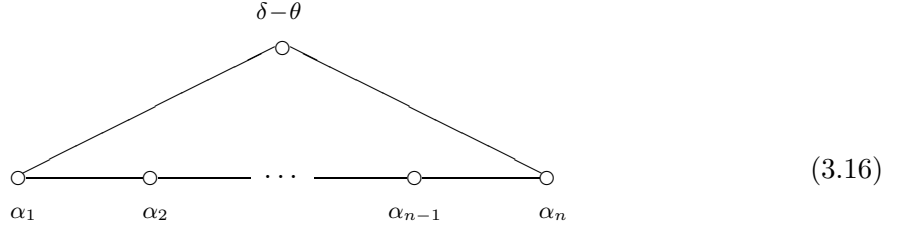
\centerline{\fs Fig.3. Dynkin diagram of the affine Lie algebra
$\wh{sl}(n\!+\!1)$.}
\\[7pt]

A general definition of the Drinfeldian $D_{q\eta}(g)$ corresponding
to a simple Lie algebra $g$ is given in \cite{T2,T3,T4}. The defining
relations for generators of $D_{q\eta}(g)$ presented in \cite{T2,T3,T4}
depend explicitly on the choice of an element
$\tl{e}_{-\th}\in U_{q}(g)$ of the weight $-\th$, such that
$g \ni \lim_{q\rightarrow 1}\tl{e}_{-\th}\neq 0$.
Here we present specification of that general definition
to the case of $g=sl(n+1)$ and set
\bn
\tl{e}_{-\th}=q^{e_{11}\!+e_{n+\!1n+\!1}}e_{n+\!11}\ .
\lb{DY18}
\ed
After some calculations we obtain the following result.
\bnpr
The Drinfeldian $D_{q\eta}'(sl(n\!+\!1))$ ($n>1$) is generated (as
a unital associative algebra over $\CC[[\log q,\eta]]$) by the algebra
$U_{q}(sl(n\!+\!1))$ and the elements $\xi_{\delta-\theta}$,
$q^{\pm h_\d}$ with the relations:
\bneqn
q^{\pm h_\d}{\rm everything}&\!\!=\!\!&{\rm everything}\,q^{\pm h_\d}~,
\lb{DY19}
\\[7pt]
q^{e_{11}}\xi_{\d-\th}&\!\!=\!\!&q^{-1}\xi_{\d-\th}q^{e_{11}}~,
\lb{DY20}
\\[7pt]
q^{e_{ii}}\xi_{\d-\th}&\!\!=\!\!&\xi_{\d-\th}q^{e_{ii}}
\qquad\quad {\rm for}\; i=2,3,\ld,n~,
\lb{DY21}
\\[7pt]
q^{e_{n\!+\!1n\!+\!1}}\xi_{\d-\th}&\!\!=\!\!&q\xi_{\d-\th}
q^{e_{n\!+\!1n\!+\!1}}~,
\lb{DY22}
\\[7pt]
[\xi_{\d-\th},e_{i\!+\!1i}]&\!\!=\!\!&0
\qquad\qquad\quad\;\;{\rm for}\; i=2,3,\ld,n-1~,
\lb{DY23}
\\[7pt]
[e_{ii\!+\!1},\xi_{\d-\th}]&\!\!=\!\!&0
\qquad\qquad\quad\;\;{\rm for}\; i=2,3,\ld,n-1~,
\lb{DY24}
\\[7pt]
[e_{12},[e_{12},\xi_{\d-\th}]_{q}]_{q}&\!\!=\!\!&0~,
\lb{DY25}
\\[7pt]
[[\xi_{\d-\th},e_{nn\!+\!1}]_{q},e_{nn\!+\!1}]_{q}&\!\!=\!\!&0~,
\lb{DY26}
\\[7pt]
[[e_{12},\xi_{\d-\th}]_{q},\xi_{\d-\th}]_{q}&\!\!=\!\!&\eta
q^{e_{11}+e_{n\!+\!1n\!+\!1}}\left(q^{-2}[e_{12},e_{n\!+\!11}]
\xi_{\d-\th}-e_{n\!+\!11}[e_{12},\xi_{\d-\th}]_{q}\right),
\lb{DY27}
\\[7pt]
[[\xi_{\d-\th},[\xi_{\d-\th},e_{nn\!+\!1}]_{q}]_{q}&\!\!=\!\!&\eta
q^{e_{11}+e_{n\!+\!1n\!+\!1}+1}
\left(q[e_{n\!+\!11},e_{nn\!+\!1}]\xi_{\d-\th}-
e_{n\!+\!11}[\xi_{\d-\th},e_{nn\!+\!1)}]_{q}\right).
\lb{DY28}
\edeqn
The Hopf structure of $D_{q\eta}'(sl(n\!+\!1))$ is defined by
the formulas (\ref{DY4})-(\ref{DY6}) for $U_{q}(sl(n\!+\!1))$
(i.e. $\D_{q\eta}(x)=\D_{q}(x)$, $S_{q\eta}(x)=S_{q}(x)$
for $(x\in U_{q}(g))$) and $\D_{q}(q^{\pm h_\d})=
q^{\pm h_\d}\ot q^{\pm h_\d}$, $S_{q}(q^{\pm h_\d})=q^{\mp h_\d}$.
The comultiplication and the antipode of $\xi_{\d-\th}$ are given by
\bn
\begin{array}{rcl}
&&\D_{q\eta}(\xi_{\d-\th}\!)=\xi_{\d-\th}\ot 1\!+\!
q^{e_{11}-e_{n\!+\!1n\!+\!1}\!-\!h_\d}\ot\xi_{\d-\th}\!
+\eta\Bigl(e_{n+11}q^{e_{n+1n+1}}\!\ot[e_{11}]
\\[5pt]
&&\quad\!\!+[\frac{h_\d}{2}\!+\!e_{n\!+\!1n\!+\!1}]q^{-\frac{h_\d}{2}}
\!\ot\!e_{n\!+\!11}q^{e_{n\!+\!1n\!+\!1}}
\!+\!\sum\limits_{i=2}^{n}\!e_{n\!+\!1i}q^{e_{n\!+\!1n\!+\!1}}\!\ot\!
e_{i1}q^{e_{ii}}\!\Big)\!\Bigl(\!q^{e_{11}}\!\ot\!q^{e_{11}}\!\Bigr),
\lb{DY29}
\end{array}
\ed
\bn
\begin{array}{lcr}
&&S_{q\eta}(\!\xi_{\d-\th}\!)\!=\!
-q^{h_{\d}-e_{11}\!+e_{n\!+\!1n\!+\!1}}\xi_{\d-\th}\!+\!\eta
[\frac{h_d}{2}\!+\!e_{11}\!+\!e_{n\!+\!1n\!+\!1}\!+\!1]
q^{\frac{h_{\d}}{2}-e_{11}\!+e_{n\!+\!1n\!+\!1}\!-1}\!
e_{n\!+\!11}
\\[5pt]
&&\qquad+\eta\sum\limits_{k=1}^{n}\!q^{-k}(q\!-\!q^{-1}\!)^{k-1}
\!\!\!\sum\limits_{n\ge i_{k}>i_{k-1}>\ld>i_{1}\ge 2}
e_{n\!+\!1i_{k}}e_{i_{k}i_{k-1}}\cd e_{i_{1}1}q^{-2e_{11}}.
\end{array}
\lb{DY30}
\ed
\lb{DYP1}
\edpr
\vskip -10pt
It is not difficult to check that the substitution
$\xi_{\d-\th}=q^{e_{11}\!+e_{n\!+\!1n\!+\!1}}e_{n\!+\!11}$
satisfies the relations (\ref{DY19})-(\ref{DY28}), i.e. there is
a simple homomorphism $D_{q\eta}(sl(n\!+\!1)) \to U_{q}(sl(n\!+\!1))$.
Moreover the both sides of the relations (\ref{DY27}) and (\ref{DY28})
are equal to zero independently. Therefore we can construct a "evaluation
representation" $\rho_{ev}$ of $D_{q\eta}(sl(n\!+\!1)$ in
$U_q(sl(n\!+\!1))\ot \CC[u]$ as follows
\bn
\begin{array}{rcccl}
\rho_{ev}(q^{h_\d})&=&1~,\qquad\quad
\rho_{ev}(\xi_{\d-\th})&=&uq^{e_{11}+e_{n\!+\!1n\!+\!1}}e_{n\!+\!11}~,
\\[7pt]
\rho_{ev}(q^{\pm h_i})&=&q^{\pm h_i}~,\qquad
\rho_{ev}(e_{\pm\a_i})&=&e_{\pm\a_i} \qquad\quad (1\le i\le n)~.
\end{array}
\lb{DY31}
\ed
We denote by $D_{q\eta}(sl(n\!+\!1))$ the Drinfeldian
$D_{q\eta}'(sl(n\!+\!1))$ with the central element $h_\d=0$.
It is obvious that
\bn
D_{q\eta=0}(sl(n\!+\!1))\simeq U_{q}(sl(n\!+\!1)[u])
\lb{DY32}
\ed
as Hopf algebras. If $q\to $1 then the limit Hopf algebra
$D_{q=1\eta}(sl(n\!+\!1))$ (and also $D_{q=1\eta}'(sl(n\!+\!1)$)
is isomorphic to the Yangian $Y_{\eta}(sl(n\!+\!1))$
$(Y_{\eta}'(sl(n\!+\!1))$ with $h_\d\neq 0$) \cite{T2}:
\bn
D_{q=1\eta}(sl(n+1))\simeq Y_{\eta}(sl(n+1))~.
\lb{DY33}
\ed
By setting $q=1$ in (\ref{DY19})-(\ref{DY30}), we obtain the
defining relations of the Yangian $Y_{\eta}'(sl(n\!+\!1))$ and its Hopf
structure in the Chevalley basis. This result is formulated as
the proposition.
\bnpr
The Yangian $Y_{\eta}'(sl(n\!+\!1)$ ($n>1$) is generated (as
an unital associative algebra over $\CC[\eta]$) by the algebra
$U(sl(n\!+\!1))$ and the elements $\xi_{\delta-\theta}$,
$h_\d$ with the relations:
\bneqn
[h_\d,{\rm everything}]&=&0~,
\lb{DY34}
\\[7pt]
[e_{11},\xi_{\d-\th}]&=&-\xi_{\d-\th}~,
\lb{DY35}
\\[7pt]
[e_{n\!+\!1n\!+\!1},\xi_{\d-\th}]&=&\xi_{\d-\th}~,
\lb{DY36}
\\[7pt]
[e_{ii},\xi_{\d-\th}]&=&0\qquad\quad {\rm for}\; i=2,3,\ld,n~,
\lb{DY37}
\\[7pt]
[\xi_{\d-\th},e_{i\!+\!1i}]&=&0\qquad\quad {\rm for}\; i=2,3,\ld,n-1~,
\lb{DY38}
\\[7pt]
[e_{ii+1},\xi_{\d-\th}]&=&0\qquad\quad {\rm for}\; i=2,3,\ld,n-1~,
\lb{DY39}
\\[7pt]
[e_{12}[e_{12},\xi_{\d-\th}]]&=&0~,
\lb{DY40}
\\[7pt]
[[\xi_{\d-\th},e_{nn\!+\!1}],e_{nn\!+\!1}]&=&0~,
\lb{DY41}
\\[7pt]
[[e_{12},\xi_{\d-\th}],\xi_{\d-\th}]\!\!&=\!\!&\eta
\Bigl([e_{12},e_{n\!+\!11}]\xi_{\d-\th}-
e_{n\!+\!11}[e_{12},\xi_{\d-\th}]\Bigr)~,
\lb{DY42}
\\[7pt]
[[\xi_{\d-\th}[\xi_{\d-\th},e_{nn\!+\!1}]]\!\!&=\!\!&\eta
\Bigl([e_{n\!+\!11},e_{nn\!+\!1}]\xi_{\d-\th}-
e_{n\!+\!11}[\xi_{\d-\th},e_{nn\!+\!1)}]\Bigr).
\lb{DY43}
\edeqn
The Hopf structure of the Yangian is trivial for
$U(sl(n+1))\oplus \CC h_\d\subset Y_{\eta}'(sl(n\!+\!1))$
(i.e. $\D_{\eta}(x)=x\ot 1\!+\!1\ot x$, $S_{\eta}(x)=-x$
for $x\in sl(n\!+\!1)\oplus\CC h_{\d}$) and it is not trivial
for the element $\xi_{\d-\th}$:
\bneqn
\D_{\eta}(\xi_{\d-\th})&\!\!=\!\!&\xi_{\d-\th}\ot 1+1\ot\xi_{\d-\th}
+\eta\,\Big(\frac{1}{2}h_{\d}\ot e_{n\!+\!11}
+\!\sum\limits_{i=1}^{n+1}e_{n\!+\!1i}\ot e_{i1}\Big)~,
\lb{DY44}
\\[5pt]
S_{\eta}(\xi_{\d-\th})&\!\!=\!\!&-\xi_{\d-\th}
+\eta\Big(\frac{1}{2}h_{\d}e_{n\!+\!11}
+\sum\limits_{i=1}^{n+1}\!\!e_{n\!+\!1i}e_{i1}\Big)~.\qquad
\lb{DY45}
\edeqn
\lb{DYP2}
\edpr
\vskip -10pt

An analog of the diagram (\ref{AH25}) for the Drinfeldian
$D_{q\eta}(sl(n\!+\!1))$ is presented by the picture:
\\[-15pt]
\bn
\mbox{\begin{picture}(100,60)
\put(-50,25){$D_{q\eta}(sl(n\!+\!1))$}
\put(-50,-25){\mbox{$Y_{\eta}(sl(n\!+\!1))$}}
\put(70,25){\mbox{$U_{q}(sl(n\!+\!1)[u])$}}
\put(70,-25){\mbox{$U(sl(n\!+\!1)[u])$}}
\put(23,29){\vector(1,0){40}}
\put(20,-21){\vector(1,0){42}}
\put(-46,16){\vector(0,-1){25}}
\put(75,16){\vector(0,-1){25}}
\put(-125,25){\mbox{$U_{q}(sl(n\!+\!1))\subset$}}
\put(-122,-25){\mbox{$U(sl(n\!+\!1))\subset$}}
\put(145,25){\mbox{$\supset U_{q}(sl(n\!+\!1))$}}
\put(142,-25){\mbox{$\supset U(sl(n\!+\!1))$}\,\,.}
\put(-70,0){\mbox{\fs$q\!\to\!1$}}
\put(78,0){\mbox{\fs$q\!\to\!1$}}
\put(30,33){\mbox{\fs$\eta\!\to\!0$}}
\put(30,-17){\mbox{\fs$\eta\!\to\!0$}}
\end{picture}}
\lb{DY46}
\ed
\\[4pt]
\begin{quote}
{\fs Fig.4. A diagram of the limit Hopf algebras
of the Drinfeldian $D_{q\eta}(sl(n\!+\!1))$\\
\phantom{Fig.1.} and their subalgebras. The arrows show passages
to the limits.}
\end{quote}
\vskip 10pt

\setcounter{equation}{0}
\section{Duality between $D_{q\eta}(sl(n\!+\!1)$ and
$\hat{H}_{q\eta}^{+}(l)$}

Let $V$ be the natural $(n\!+\!1)$-dimensional representation of the
quantum algebra $U_{q}(sl(n\!+\!1))$ with basis
$\{v_1,v_2,\ld,v_{n\!+\!1}\}$ on which the action of
$U_{q}(sl(n\!+\!1))$
is given by
\bn
\begin{array}{rcl}
e_{i\!-\!1i}v_k&=&\d_{ik}v_{k\!+\!1}~,
\\[7pt]
e_{i\!+\!1i}v_k&=&\d_{ik}v_{k\!-\!1}~,
\\[7pt]
q^{\pm e_{ii}}v_k&=&q^{\pm\d_{ik}}v_{k}~.
\end{array}
\lb{D1}
\ed
Let $T$: $V\ot V\to V\ot V$ be a linear map given by
\bn
T(v_r\ot v_s)=\left\{\begin{array}{rcll}
qv_r\!\!\!&\ot\!\!\!&v_s& {\rm if}\quad r=s~,\\
v_s\!\!\!&\ot\!\!\!&v_r&{\rm if}\quad r\le s~,\\
v_s\!\!\!&\ot\!\!\!&v_r+(q-q^{-1})v_r\ot v_s&{\rm if}\quad r\ge s~.
\end{array}\right.
\lb{D2}
\ed
It is not difficult to check that the elements
$\sg_i\in {\rm End}_{\CC}(V^{\ot l})$
which act as $T$ on $i^{-th}$ and $(i\!+\!1)^{-th}$ factors of the tensor
product, and as the identity on the other factors, for $i=1,2,\ld, l$
define the representation of the Hecke algebra $H_{q}(l)$ on $V^{\ot l}$.

We say that a representation of $D_{q\eta}(sl(n\!+\!1)$ has a level $l$
if its restriction to $U_{q}(sl(n\!+\!1))$ is sum of representations
each of which occurs in $V^{\ot l}$. Now we announce the main result.
\bnth
(i) Let $M$ be a finite-dimensional right $\hat{H}_{q\eta}^{+}(l)$-module
and we set $W_M=M\ot_{H_{q}(l)}V^{\ot l}$. Then there exists a
homomorphism $\pi$: $D_{q\eta}(sl(n\!+\!1)\to {\rm End}_{\CC} W_M$
such that
\bneqn
\pi(x)(m\ot{\bf v})&=&m\ot\D_{q}^{(l)}(x){\bf v}\qquad\quad
{\rm for}\; x\in U_{q}(sl(n\!+\!1))~,
\lb{D3}
\\[7pt]
\pi(\xi_{\d-\th})(m\ot{\bf v})&=&
m\ot\Big(\D_{q\eta}^{(l)}
(\xi_{\d-\th})\Bigr|_{\stackrel{i}{\xi}_{\d-\th}=u_i}\Big){\bf v}
\lb{D4}
\edeqn
for $m\in M,\;{\bf v}\in V^{\ot l}$. For $l\le n$ the functor
$\cF_{q\eta}(M)$: $M\to W_M$ is an equivalence between the category
of finite-dimensional right $\hat{H}_{q\eta}^{+}(l)$-modules
and the category of  finite-dimensional left
$D_{q\eta}(sl(n\!+\!1))$-modules of level $l$.

\nin
(ii) For $\eta=0$ the functor $\cF_{q\eta=0}(M)$ is an equivalence
between the category of finite-dimensional right
$\hat{H}_{q}^{+}(l)$-modules and the category of finite-dimensional
left $U_{q}(sl(n\!+\!1)$-modules of level $l\le n$.

\nin
(iii) For $q\to 0$ the functor $\cF_{q=\!1\eta}(M)$ is
an equivalence between the category of finite-dimensional right
$\Lm_{\eta}(l)$-modules and the category of
finite-dimensional left $Y_{\eta}(sl(n\!+\!1)$-modules of level $l\le n$.
\lb{DT1}
\edth
\vskip -10pt
\nin
Here $\D_{q\eta}^{(l)}$ is the $l$-fold coproduct
\bn
\D_{q\eta}^{(l)}:D_{q\eta}(sl(n\!+\!1))\to D_{q\eta}(sl(n\!+\!1))\ot\cd
\ot D_{q\eta}(sl(n\!+\!1))\qquad (l-{\rm fold})~.
\lb{D5}
\ed
In particular
\bn
\D_{q\eta}^{(2)}(\cdot)=\D_{q\eta}(\cdot)
\lb{D6}
\ed
The symbol $\stackrel{i}{\xi}_{\d-\th}=u_i$ in (\ref{D4}) means that
the $i^{-th}$ component of the affine element $\xi_{\d-\th}$ in the
$l$-fold coproduct $\D_{q\eta}^{(l)}(\xi_{\d-\th})$ has to replace
by the affine Hecke element $u_i$.

The proof of the part $(i)$ of Theorem \ref{DT1} is analogous to
the proof of the duality theorem between the affine Hecke algebra
$\hat{H}_{q\eta}(l)$ and the quantum affine algebra
$U_{q}(\wh{sl}(n\!+\!1))$ (see \cite{ChP1}).
The parts $(ii)$ and $(iii)$ are proven by direct
comparison of $\cF_{q\eta=0}(M)$ and $\cF_{q=\!1\eta}(M)$ with
the Chari-Pressley's and Drinfeld's functors {\cite{ChP1,D2}}.

\section*{Acknowledgments}

The first author (V.N.T.) is grateful to
the Org. Committee of the Intern. Symposium
"Quantum Theory and Symmetries",
H.-D. Doebner, V.K. Dobrev, J.-D. Hennig and W. L\"ucke,
for the support of his visit on this Workshop
and he is also thankful to Toulon University for the support
of his visit to the Marseille Center of Theoretical Physics CNRS
where the first stages of the work were done.
This work was supported by the program of French--Russian scientific
cooperation (CNRS grant PICS-608 and grant RFBR-98-01-22033)
and also by grant RFBR-98-01-00303 (V.N. Tolstoy).


\begin{thebibliography}{99}

\bibitem[1]{BZ}
Bernshtein I.N. and Zelevinskii A.V.,
Representation of the group $GL(n,F)$, where $F$ is local
Archimedean field, {\it Usp.Math. Nauk,} {\bf 31} (1976), 5-70.

\bibitem[2]{ChP1}
Chari V. and Pressley A.,
Quantum Affine Algebras and Affine Hecke Algebras,
{\it e-print math.QA/9501003} (1995).

\bibitem[3]{ChP2}
Chari V. and Pressley A.,
{\it A guide to quantum groups}, Cambridge University Press, 1994.

\bibitem[4]{D1}
Drinfeld V.G.,
Hopf algebras and quantum Yang-Baxter equation,
{\it Soviet Math. Dokl.,} {\bf 283} (1985), 1060-1064.

\bibitem[5]{D2}
Drinfeld V.G.,
Degenerate affine Hecke algebras and Yangians,
{\it Func. Anal. Appl.,} {\bf 20} (1986), 62-64.

\bibitem[6]{J1}
Jimbo M.,
A q-difference analogue of $U(g)$ and Yang-Baxter equation,
{\it Lett. Math. Phys.,} {\bf 10} (1985), 63-69.

\bibitem[7]{J2}
Jimbo M.,
A q-analogue of $U_{q}(gl(n\!+\!1))$, Hecke algebras and Yang-Baxter
equation, 
{\it Lett. Math. Phys.,} {\bf 11} (1986), 247-252.

\bibitem[8]{KT1}
Khoroshkin S.M. and Tolstoy V.N.,
Universal R-matrix for quantized (super)algebras,
{\it Commun. Math. Phys.,}{\bf 141} (1991), 599-617.

\bibitem[9]{KT2}
Khoroshkin S.M., and Tolstoy  V.N.,
Twisting of quantum (super)algebras. Connection of Drinfeld's
and Cartan-Weyl realizations for quantum affine algebras,
{\it Preprint MPIM Bonn (Germany), MPI/94-23} (1994), 29p.;
{\it e-print hep-th/9404036} (1994).

\bibitem[10]{R}
Rogavski J.D.,
On modules over the Hecke algebras of p-adic group,
{\it Invent. Math.,} {\bf 79} (1985), 443-465.

\bibitem[11]{T1}
Tolstoy V.N.,
Extremal projectors for quantized Kac-Moody superalgebras and some
of their applications,
{\it Lectures Notes in Physics 370} (1990), 118-125.

\bibitem[12]{T2}
Tolstoy V.N.,
Connection between Yangians and Quantum Affine Algebras,
Proceedings of the X-th Max Born Symposium,
(Wroclav, 1996, eds: J. Lukierski, M.Mozrzymas).
PWN - Polish Scientific Publishers - Warszawa (1997), 99-117.

\bibitem[13]{T3}
Tolstoy V.N.,
Two-parameter deformations of loop algebras and superalgebras,
Proc. of the 5-th Wigner Symposium,
(Vienna, 1997, eds: P.Kasperkovitz, D.Grau).
World Scientific, Singapore-New Jersey-London-Hong Kong
(1998), 25-27; 
{\it e-print math.QA/9712028} (1997).

\bibitem[14]{T4}
Tolstoy V.N.,
Drinfeldians,
Proc. vol. {\it "Lie Theory and its Application in Physics II"}
(eds: H.-D. Doebner, V.K. Dobrev and J. Hilgert, World Scientific,
Singapore, 1998, 981-02-3539-9, pp. 325-337;
{\it e-print math.QA/9803008} (1998).

\end{thebibliography}
\end{document}